\documentclass[10pt, reqno]{amsart}
\usepackage{amsmath, amssymb}
\usepackage{dsfont}
\usepackage{amscd}
\usepackage[mathscr]{euscript} 
\usepackage[all]{xy}
\usepackage{url}
\usepackage{stmaryrd}
\usepackage{comment}
\usepackage[retainorgcmds]{IEEEtrantools}

\usepackage{xcolor}
\usepackage[colorlinks=true]{hyperref}

\usepackage{graphicx}
\usepackage{mathtools}
\usepackage[dvipsnames]{xcolor}
\usepackage{centernot}
\usepackage{marvosym}
\usepackage{tikz}
\usepackage{tikz-cd}
\usetikzlibrary{calc, arrows.meta, positioning}
\usepackage{xspace}
\usepackage{bbm}
\usepackage[title]{appendix}

\usepackage{tkz-fct}
\usepackage{soul}
\usepackage{lineno}
\usepackage[T1]{fontenc}
\usepackage[utf8]{inputenc}

\usepackage[english]{babel}
\usepackage[autostyle]{csquotes}
\usepackage[shortlabels]{enumitem}

\usepackage[backend=bibtex,  sorting=nyt, style=alphabetic]{biblatex}
\allowdisplaybreaks 

\title[IMC and solvable groups]
{The Idempotent Measure Conjecture beyond solvable and stable groups}


\usepackage{orcidlink}

\author[D. M. HOFFMANN]{Daniel Max Hoffmann$^{\dagger}$\
\orcidlink{0000-0002-4514-269X}}
\thanks{$^{\dagger}$SDG. The author is supported
by the National Science Centre (Narodowe Centrum Nauki, Poland) 
grant no. 2021/43/B/ST1/00405,
and by the UW IDUB PSP no. 501-D110-20-0001381.}
\address{$^{\dagger}$
Instytut Matematyki\\
Uniwersytet Warszawski\\
Warszawa\\
Poland}
\email{daniel.max.hoffmann@gmail.com}
\urladdr{{https://sites.google.com/site/danielmaxhoffmann/}}

\date{\today}

 \DeclareMathOperator{\aut}{Aut} 

\DeclareMathOperator{\stab}{Stab}

\DeclareMathOperator{\tp}{tp}

\DeclareMathOperator{\inv}{inv}

\newcommand{\Aut}{\operatorname{Aut}}

\DeclareMathOperator{\eq}{eq}

\newtheorem{theorem}{Theorem}[section]

\newtheorem{lemma}[theorem]{Lemma}
\newtheorem{cor}[theorem]{Corollary}
\newtheorem{fact}[theorem]{Fact}

\theoremstyle{definition}
\newtheorem{definition}[theorem]{Definition}

\theoremstyle{remark}
\newtheorem*{theorem*}{Theorem}
\newtheorem*{cor*}{Corollary}

\theoremstyle{definition}
\theoremstyle{definition}

\theoremstyle{definition}

\theoremstyle{remark}

\newtheorem*{clm*}{Claim}

\AtEndEnvironment{proof}{\setcounter{clm}{0}}

\makeatletter
\providecommand*{\cupdot}{%
  \mathbin{%
    \mathpalette\@cupdot{}%
  }%
}
\newcommand*{\@cupdot}[2]{%
  \ooalign{%
    $\m@th#1\cup$\cr
    \sbox0{$#1\cup$}%
    \dimen@=\ht0 %
    \sbox0{$\m@th#1\cdot$}%
    \advance\dimen@ by -\ht0 %
    \dimen@=.5\dimen@
    \hidewidth\raise\dimen@\box0\hidewidth
  }%
}

\providecommand*{\bigcupdot}{%
  \mathop{%
    \vphantom{\bigcup}%
    \mathpalette\@bigcupdot{}%
  }%
}
\newcommand*{\@bigcupdot}[2]{%
  \ooalign{%
    $\m@th#1\bigcup$\cr
    \sbox0{$#1\bigcup$}%
    \dimen@=\ht0 %
    \advance\dimen@ by -\dp0 %
    \sbox0{\scalebox{2}{$\m@th#1\cdot$}}%
    \advance\dimen@ by -\ht0 %
    \dimen@=.5\dimen@
    \hidewidth\raise\dimen@\box0\hidewidth
  }%
}
\makeatother

\def\Ind#1#2{#1\setbox0=\hbox{$#1x$}\kern\wd0\hbox to 0pt{\hss$#1\mid$\hss}
\lower.9\ht0\hbox to 0pt{\hss$#1\smile$\hss}\kern\wd0}

\def\notind#1#2{#1\setbox0=\hbox{$#1x$}\kern\wd0
\hbox to 0pt{\mathchardef\nn=12854\hss$#1\nn$\kern1.4\wd0\hss}
\hbox to 0pt{\hss$#1\mid$\hss}\lower.9\ht0 \hbox to 0pt{\hss$#1\smile$\hss}\kern\wd0}

\newcommand{\FC}{\mathfrak{C}}

\newcommand{\CL}{{\mathcal L}}

\newcommand{\Av}{\mathrm{Av}}
\newcommand{\supp}{\mathrm{supp}}

\begin{document}

\begin{abstract}
Pym's theorem \cite{Pym62} states that a regular Borel probability
measure on a locally compact group is idempotent if and only if
it is the normalized Haar measure of a compact subgroup.
The Idempotent Measure Conjecture, the definable-group counterpart
of \cite[Conjecture~(A)]{GHK}, proposes an analogous characterization:
a global fim Keisler measure on a definable group is idempotent
if and only if it is the unique global Keisler measure concentrating
on a type-definable subgroup and invariant under translation
by elements of that subgroup.

We prove this conjecture for virtually solvable definable groups
in arbitrary theories.
More generally, we establish it for definable groups admitting
a finite series of definable normal subgroups whose successive
quotients are either virtually solvable or stably embedded
with stable induced structure.

Finally, we disprove \cite[Conjecture~(A)]{GHK}
in its automorphism-group
formulation in the theory of dense linear
orders without endpoints.
\end{abstract}

\maketitle

\section{Introduction}
The classification of idempotent measures is a central problem in the
study of convolution semigroups 
\cite{Wendel,Rudin,Glicksberg2,Pym62}. 
In the classical setting, a regular
Borel probability measure on a locally compact group satisfies
$\mu*\mu=\mu$, i.e. $\mu$ is idempotent,
if and only if it is the normalized Haar measure of a
compact subgroup \cite[Theorem~A.4.1]{Pym62}. 
The Idempotent Measure Conjecture asks for a
corresponding characterization in model theory, with type-definable
subgroups replacing compact subgroups.

Let $G$ be a group definable in a monster model $\FC$ of a complete
theory $T$. A global Keisler measure on $G$ is a finitely additive
probability measure on its definable subsets (cf. \cite{KyleThesis} for a nice exposition of the subject). Definable convolution,
studied in \cite{Artem_Kyle,Artem_Kyle2,CGK}, is obtained by pushing
the Morley product of measures forward under group multiplication.
We consider frequency interpretation measures, or fim measures,
introduced in \cite[Definition~2.7]{HruPiSi13}.
In NIP theories, these are precisely the generically stable measures
of \cite{Anand_Udi2011}; see \cite[Theorem~3.2]{HruPiSi13}.
Generic stability for measures beyond NIP has been studied in
\cite{ConGann20,InTheWild,ConGanHan23,CGH26}.
Related developments concerning extension domination and the
compatibility of randomization transfer maps with Morley products
appear in \cite{GannonYeDomination,GannonTransfer}.
The precise definition of a fim measure 
is recalled in Section~\ref{sec:preliminaries}.
The Idempotent Measure Conjecture
predicts that a global fim measure is idempotent if
and only if it is the unique global Keisler measure concentrating on
a type-definable subgroup and invariant under right translation by
elements of that subgroup; see \cite[Conjecture~1.3]{GHK}.

For a global fim measure $\mu$, its right stabilizer
$$\stab_{\mathrm r}(\mu) =\{g\in G:\mu*\delta_g=\mu\}$$
is type-definable (cf. \cite[Proposition~5.3]{Artem_Kyle}, \cite[Fact~3.35]{CGK}). For an idempotent such measure, the conjectural
characterization is equivalent to concentration on this stabilizer,
a property called \emph{generic transitivity}
\cite[Proposition~3.37 and Remark~3.39]{CGK}.
Accordingly, the conjecture is formulated in
\cite[Problem~3.41]{CGK} as the question whether every global fim
idempotent measure is generically transitive. 
Chernikov and Gannon characterized idempotent measures in stable
theories as translation-invariant measures on type-definable
subgroups \cite[Theorem~5.8]{Artem_Kyle},
and Chernikov, Gannon and Krupi\'nski
established it for abelian groups in arbitrary theories
\cite[Theorem~3.45]{CGK}.

We prove the conjecture for virtually solvable definable groups in
arbitrary theories (Corollary~\ref{cor:virtually-solvable}). Here
virtual solvability is an abstract group-theoretic assumption: the
solvable subgroup of finite index need not be definable. We also
combine this result with the stable case along finite normal series.
Our main theorem is the following; see Theorem~\ref{thm:main}.

\begin{theorem}\label{thm:intro-main}
Let $G$ be a definable group admitting a finite series
$$\{1\}=G_0\trianglelefteq G_1\trianglelefteq\cdots \trianglelefteq G_m=G$$
of definable subgroups normal in $G$. Suppose that every successive
quotient $G_{i+1}/G_i$ is either virtually solvable or stably embedded
with stable induced structure. Then every global fim idempotent
Keisler measure on $G$ is generically transitive.
\end{theorem}

The proof proceeds by lifting generic transitivity from quotient
measures. Its main ingredient is the following criterion
(Theorem~\ref{thm:commuting-kernels}).

\begin{theorem}\label{thm:intro-commuting-kernels}
Let $N_1,N_2\trianglelefteq G$ be definable normal subgroups such
that $[N_1,N_2]=1$, and let $h_i:G\to G/N_i$ be the quotient maps.
If $\mu$ is a global fim idempotent measure on $G$ and both
$(h_1)_*\mu$ and $(h_2)_*\mu$ are generically transitive, then
$\mu$ is generically transitive.
\end{theorem}

Taking $N_1=N$ and $N_2=C_G(N)$ gives a criterion involving the
quotients $G/N$ and $G/C_G(N)$. If $N$ is abelian, the latter
quotient is a homomorphic image of the former, so generic
transitivity lifts through abelian kernels. An induction on derived
length, using definable characteristic abelian subgroups, gives the
same conclusion for solvable kernels. Finite kernels are treated
using the finiteness of $G/C_G(N)$. Passing to the \emph{definable solvable
radical} (Lemma \ref{lem:solvable-radical}) then yields transfer through virtually solvable kernels
(Theorem~\ref{thm:virtually-solvable-kernel-transfer}).

For a stably embedded normal subgroup $N$ with stable induced
structure, the conjugation action identifies $G/C_G(N)$ with a group
definable in the imaginary sorts of that induced structure.
Comparing ambient and induced Morley products and convolutions
allows us to apply the stable classification of \cite{Artem_Kyle} to this quotient.
The same criterion therefore gives transfer through such kernels
(Theorem~\ref{thm:stable-kernel-transfer}). Iterating the two transfer
theorems proves Theorem~\ref{thm:intro-main}.

There is also an automorphism version of the conjecture.
\cite{GHK} introduced convolution on $M$-invariant global Keisler
measures extending $\tp(\bar{m}/\emptyset)$ for an enumeration $\bar{m}$ of a fixed small model $M$,
and formulated the analogous
characterization using relatively type-definable subgroups of the
automorphism group \cite[Conjecture~(A)]{GHK}. The affine sort
construction relates this operation to definable convolution.
Subsequently, a conditional characterization in this setting was obtained under left-simplicity of the support and a local invariance condition \cite[Theorem~3.14]{HRz}.
In Appendix~\ref{sec:automorphism}, we show that Conjecture~(A)
fails in general: in the theory of dense linear orders without
endpoints, we construct a smooth fim idempotent measure for which
the stabilizer has measure zero (Theorem~\ref{thm:auto}). 

Section~\ref{sec:preliminaries} collects the measure-theoretic
preliminaries. Section~\ref{sec:quotient-transfer} proves the
criterion for commuting kernels. Section~\ref{sec:solvable}
establishes the results for virtually solvable groups, and
Section~\ref{sec:normal-series} proves the stable-kernel transfer
theorem and the main theorem. Appendix~\ref{sec:automorphism} is dedicated to the counterexample for automorphism variant of \cite[Conjecture (A)]{GHK}.

\ \\ \\
\textbf{Use of AI.} The author used ChatGPT (GPT-6 Astra)
to assist in identifying examples,  developing proofs, and drafting the text. 
In particular, 
the exact shape of Lemma \ref{lem:stable-coding} and the counterexample in Appendix were suggested by AI.
The author subsequently refined and streamlined the arguments and assume full responsibility for the final content.

\section{Preliminaries}\label{sec:preliminaries}
We work in a complete theory $T$ with a monster model $\FC$.
Sets of parameters and elementary submodels are small unless otherwise
specified. We work in $T^{\eq}$ when taking definable quotients.
All groups and maps are definable, with parameters allowed, unless
explicitly stated otherwise. Every named base contains the parameters
defining the groups and maps under consideration. 

We identify each definable group $G$ with its group of
$\FC$-points. If $G$ is defined over a small model $M$,
then $G(M)$ denotes its group of $M$-points.
For a definable group $G$, its identity is
denoted by $1$. We use the convention $[g,h]=g^{-1}h^{-1}gh$ and
$g^h=h^{-1}gh$.

\subsection{Measures and Morley products}
For a tuple of variables $x$ and a parameter set $A$, let
$\mathfrak M_x(A)$ be the space of Keisler measures over $A$ in the
variables $x$. We identify such a measure with the corresponding
regular Borel probability measure on $S_x(A)$, and write
$\supp(\mu)$ for its support. A measure is \emph{global} if its domain
is the Boolean algebra of formulas over $\FC$. We write
$\mathfrak M_G(A)$ for measures concentrated on $G$.
If $P$ is type-definable, saying that $\mu$ \emph{concentrates on} $P$
means that $\mu(X)=1$ for every definable set $X$ containing $P$.
For a Keisler measure $\mu$, we write ``$\mu_x$'' to indicate that it is considered as a measure in the variables $x$.

A global measure $\mu$ is \emph{$M$-invariant} if
$\mu(\varphi(x;b))$ depends only on $\tp(b/M)$, for every formula
$\varphi(x;y)$. It is \emph{definable over $M$} if, in addition, each
function
$$F^{\varphi}_{\mu,M}:S_y(M)\longrightarrow[0,1],\qquad \tp(b/M)\longmapsto\mu(\varphi(x;b)),$$
is continuous. These functions may be uniformly approximated by
finite-valued $M$-definable functions.

For a global definable measure $\mu$ and an arbitrary global measure
$\nu$, the \emph{Morley product} is given by
$$(\mu_x\otimes\nu_y)(\varphi(x,y;c)) =\int_{S_y(M)}\mu(\varphi(x;b;c))\,d(\nu|M)(b),$$
where $M$ contains $c$ and is a model over which $\mu$ is definable.
The value is independent of the choice of $M$. Thus the left factor is
integrated first, as in \cite[Definition~2.23]{KyleThesis}.
We set $\mu^{(1)}=\mu$ and
$$\mu^{(n+1)}_{x_1,\ldots,x_{n+1}} =\mu_{x_{n+1}}\otimes\mu^{(n)}_{x_1,\ldots,x_n}.$$
In particular, $\mu^{(n)}$ has $n$ coordinates. The compatible union of
these measures is denoted by $\mu^{(\omega)}$.

\begin{definition}\label{def:fim}
A global measure $\mu_x$, definable over $M$, is a
\emph{frequency interpretation measure over $M$}, or \emph{fim over $M$},
if, for every
formula $\varphi(x;y)$, there are formulas
$\theta_n(x_1,\ldots,x_n)\in\CL(M)$ such that
$\mu^{(n)}(\theta_n)\to1$ and, for every $\varepsilon>0$, for all
sufficiently large $n$, every tuple $\bar{a}=(a_1,\ldots,a_n)$ satisfying $\theta_n$ satisfies
$$\sup_b\left| \Av(\bar{a})(\varphi(x;b)) -\mu(\varphi(x;b))\right|<\varepsilon.$$
Here $\Av(\bar a)(\varphi(x;b))=\frac{1}{n} |\{i\in\{1,\ldots,n\}:\FC\models\varphi(a_i;b)\}|$. A measure is \emph{fim} if it is fim over some small model.
\end{definition}

This is the usual notion of fim, using the definability of fim measures;
see \cite[Definition~3.5]{CGK}.
A measure fim over $M$ remains fim over
every elementary submodel containing $M$. We will enlarge the base
model to include the parameters of the maps and type-definable sets
under consideration.

If $f$ is definable, its pushforward is
$f_*\mu(X)=\mu(f^{-1}(X))$. We use the following standard properties.

\begin{fact}\label{fact:functorial}
Morley products of definable measures are definable and associative.
The Morley product of two global measures fim over $M$ is fim over $M$.
Closure of fim measures under Morley products is established in
\cite[Theorem~1.3]{CGH26}.
A fim measure commutes, under the Morley product, with every global
definable measure. Definable pushforwards preserve definability and fim,
and for definable maps $f$ and $g$,
$$(f\times g)_*(\mu\otimes\nu)=(f_*\mu)\otimes(g_*\nu)$$
whenever the displayed products are defined with $\mu$ definable.
Consequently, the Morley powers of a fim measure are symmetric under
permutations of their coordinates.
\end{fact}

For commutation and pushforwards, see
\cite[Fact~3.7 and Proposition~3.26]{CGK}; associativity for definable
measures is established in \cite[Theorem~2.18]{InTheWild};
see also \cite[Proposition~2.24]{KyleThesis} for definability of products.
The displayed functoriality identity follows directly by applying the
defining integral to a formula on the target. 

\begin{fact}[\cite{CGK}, Proposition~3.11(1)]\label{fact:uniqueness}
Suppose that $\mu$ is fim over $M$, that $B\supseteq M$, and that
$\nu$ is a global $B$-invariant definable measure. If
$\nu|B=\mu|B$, then $\nu=\mu$.
\end{fact}

\begin{lemma}\label{lem:parameter-pushforward}
Let $\mu_x$ be a global measure definable over $M$, let $\lambda_z$
be fim over $M$, and let $\nu_y$ be a global measure. Suppose that
$f(x,y)$ is $M$-definable and
$$((x,y)\mapsto(f(x,y),y))_*(\mu\otimes\nu) =\lambda\otimes\nu.$$
Then $(f_b)_*\mu=\lambda$ for every
$b\models p\in\supp(\nu|M)$, where $f_b(x)=f(x,b)$.
\end{lemma}
\begin{proof}
Fix $\varphi(z,y)\in\CL(M)$. The continuous functions on $S_y(M)$
given by
$$b\longmapsto\mu_x(\varphi(f(x,b),b)),\qquad b\longmapsto\lambda_z(\varphi(z,b))$$
have equal integrals over every clopen set, by the assumed identity.
They therefore agree on $\supp(\nu|M)$.

Fix $b$ with $\tp(b/M)\in\supp(\nu|M)$, and let
$\psi(z)\in\CL(Mb)$. We may write $\psi(z)=\varphi(z,b)$ for some
$\varphi(z,y)\in\CL(M)$. The preceding equality at $\tp(b/M)$ gives
$$((f_b)_*\mu)(\psi(z))=\mu_x\bigl(\varphi(f(x,b),b)\bigr) =\lambda_z\bigl(\varphi(z,b)\bigr)=\lambda(\psi(z)).$$
Since $\psi$ was arbitrary, $(f_b)_*\mu$ and $\lambda$ agree over
$Mb$. The former is definable over $Mb$, so
Fact~\ref{fact:uniqueness} gives the result.
\end{proof}

The following lemma is a consequence of
\cite[Theorem~3.13]{CGK}. We give a direct proof from the
definition of fim.

\begin{lemma}\label{lem:finite-couplings}
Let $\mu_x$ be fim over $M$, let $\nu_y$ be a global measure,
and let $\tau_{x,y}$ be a global measure. Suppose that, for arbitrarily
large $n$, there is a global measure
$\lambda_n(x_1,\ldots,x_n,y)$ such that
$$\lambda_n|_{x_1,\ldots,x_n}=\mu^{(n)},\qquad \lambda_n|_y=\nu,\qquad \lambda_n|_{x_i,y}=\tau\quad(1\leq i\leq n).$$
Then $\tau=\mu\otimes\nu$.
\end{lemma}

\begin{proof}
Fix a formula $\varphi(x,y;c)$ and $\varepsilon>0$.
Choose a small model $N\supseteq Mc$, and let
$F:S_y(N)\to[0,1]$ be the continuous function given by
$$F(\tp(b/N))=\mu(\varphi(x,b;c)).$$
Applying the fim property to $\varphi(x;y,z)$, we obtain, for all
sufficiently large $n$, a formula
$\theta_n(x_1,\ldots,x_n)\in\CL(M)$ such that
$\mu^{(n)}(\theta_n)>1-\varepsilon$ and
$$\left| \Av(a_1,\ldots,a_n)(\varphi(x,b;c)) -\mu(\varphi(x,b;c)) \right|<\varepsilon$$
whenever $\theta_n(a_1,\ldots,a_n)$ holds, for every $b$.

Choose such an $n$ for which $\lambda_n$ exists as in the hypothesis.
Let
$$\pi_y:S_{x_1,\ldots,x_n,y}(N)\longrightarrow S_y(N)$$
be the restriction map, and define
$$A_n:=\frac1n\sum_{i=1}^n \mathbf 1_{[\varphi(x_i,y;c)]}.$$
The preceding estimate gives $|A_n-F\circ\pi_y|<\varepsilon$
on $[\theta_n]$, while $|A_n-F\circ\pi_y|\leq1$ everywhere.
Using the prescribed restrictions of $\lambda_n$ and the definition
of the Morley product, we obtain
$$
\begin{aligned}
\left| \tau(\varphi(x,y;c)) -(\mu\otimes\nu)(\varphi(x,y;c)) \right| &= \left| \int_{S_{x_1,\ldots,x_n,y}(N)} (A_n-F\circ\pi_y)\,d(\lambda_n|N) \right|\\
&\leq \varepsilon+\lambda_n(\neg\theta_n) =\varepsilon+1-\mu^{(n)}(\theta_n) <2\varepsilon.
\end{aligned}
$$
Since $\varepsilon>0$ and $\varphi(x,y;c)$ were arbitrary,
$\tau=\mu\otimes\nu$.
\end{proof}

\subsection{Convolution and generic transitivity}
For global definable measures on $G$, convolution is the pushforward
of the Morley product under multiplication:
$$\mu*\nu=m_*(\mu\otimes\nu),\qquad m(x,y)=xy.$$
This is the convention of \cite[Definition~6.9]{KyleThesis}.
The Dirac measure at $g$ is denoted by $\delta_g$. Thus
$\mu*\delta_g$ and $\delta_g*\mu$ are the right and left translates,
respectively. We write $\eta^{*n}$ for the $n$-fold convolution power,
and set $\eta^{*0}=\delta_1$.
Convolution is associative on global definable measures, by
associativity of the Morley product
\cite[Theorem~2.18]{InTheWild}, functoriality
(Fact~\ref{fact:functorial}), and associativity of group multiplication.
Moreover, if $f:G\to H$ is a definable homomorphism, then
$$f_*(\mu*\nu)=(f_*\mu)*(f_*\nu).$$
Indeed, writing $m_G$ and $m_H$ for the respective multiplication
maps, we have $f\circ m_G=m_H\circ(f\times f)$, so
Fact~\ref{fact:functorial} gives
$$
\begin{aligned}
f_*(\mu*\nu) &= (f\circ m_G)_*(\mu\otimes\nu)\\
&= (m_H)_*\bigl((f\times f)_*(\mu\otimes\nu)\bigr)\\
&= (m_H)_*\bigl((f_*\mu)\otimes(f_*\nu)\bigr)\\
&= (f_*\mu)*(f_*\nu).
\end{aligned}
$$

The \emph{right stabilizer} of a global measure $\mu$ is
$$\stab_{\mathrm r}(\mu)=\{g\in G(\FC):\mu*\delta_g=\mu\}.$$
It is a subgroup and, if $\mu$ is definable over $M$, it is
$M$-type-definable \cite[Fact~3.35]{CGK}.

\begin{definition}\label{def:gt}
A global fim idempotent $\mu$ on $G$ is \emph{generically transitive}
if it concentrates on $\stab_{\mathrm r}(\mu)$.
\end{definition}

\begin{fact}\label{fact:gt}
Let $\mu$ be a global fim idempotent on $G$, and put
$H=\stab_{\mathrm r}(\mu)$.
\begin{enumerate}[(i)]
\item Generic transitivity is equivalent to
$$((x,y)\mapsto(xy,y))_*(\mu^{(2)})=\mu^{(2)}.$$
\item If $\mu$ is generically transitive, it is the unique left-invariant
and the unique right-invariant global measure concentrating on $H$.
In particular, inversion preserves $\mu$.
\item Definable homomorphic pushforwards of generically transitive fim
idempotents are generically transitive.
\item If $\mu$ is generically transitive, then
$$\{g:\mu*\delta_g*\mu=\mu\}=H.$$
\end{enumerate}
\end{fact}
\begin{proof}
The first two assertions are
\cite[Proposition~3.37 and Remark~3.39]{CGK}; inversion invariance
follows from the uniqueness of the invariant measures on $H$.
The third follows from the first and Fact~\ref{fact:functorial}.
For the last assertion, the inclusion from right to left is immediate.
Suppose that $g\notin H$. 
The product measure $\mu\otimes\mu$ concentrates on $H\times H$.
Consequently,
$$\mu*\delta_g*\mu =((x,y)\mapsto xgy)_*(\mu\otimes\mu)$$
concentrates on $HgH$. The latter set is type-definable, since it
is the image of the type-definable set $H\times H$ under a definable
map.
Since $H$ and $HgH$ are disjoint type-definable
sets, $\mu$ concentrates on $H$ and $\mu*\delta_g*\mu$ concentrates
on $HgH$, we obtain $\mu*\delta_g*\mu\neq\mu$.
\end{proof}

\begin{fact}\label{fact:finite}
Every idempotent probability measure on a finite group is the
normalized counting measure on a subgroup. In particular, it is fim
and generically transitive.
\end{fact}

This is the finite case of the classical idempotent measure theorem
\cite{kawada1940probability}. We will also use its stable model-theoretic counterpart.

\begin{fact}[\cite{Artem_Kyle}, Theorem~5.8]\label{fact:stable-imc}
In a stable theory, every global idempotent measure on a definable
group is the unique invariant measure concentrating on a
type-definable subgroup. In particular, every global fim idempotent
on such a group is generically transitive.
\end{fact}

\subsection{Induced structures}
Fix a small parameter set $B$ and a $B$-definable set $D$.
By the \emph{induced structure on $D$ over $B$} we mean the structure
with a relation for every $B$-definable subset of every finite power
of $D$. This is a small language. We say that $D$ is
\emph{stably embedded over $B$} if every ambient-definable subset of
every finite power of $D$ is definable, with parameters from $D$, in
this induced structure \cite{PillayStableEmbeddedness}. When we refer to a stably embedded set with
stable induced structure, these conditions are understood over some
fixed small base.

\section{Generic transitivity under quotients}\label{sec:quotient-transfer}

We first establish the identities needed to compare translates of an
idempotent measure on the fibres of a definable quotient.

The first assertion of the following lemma is the measure analogue
of the sandwich identity for invariant idempotent types in
\cite[Remark~2.22]{CGK}.

\begin{lemma}\label{lem:sandwich-locus}
Let $G$ be an $M$-definable group, and let 
$\mu\in\mathfrak M_G(\FC)$ be fim over $M$ and idempotent. 
Then
$$((x,y,z)\mapsto(xyz,y))_*(\mu^{(3)})=\mu^{(2)}.$$
Moreover, the set
$$K_\mu=\{b\in G:\mu*\delta_b*\mu=\mu\}$$
is $M$-type-definable. For every
$b\models p\in\supp(\mu|M)$, one has $b\in K_\mu$, and both
$\mu$ and $\mu*\delta_b$ concentrate on $K_\mu$.
\end{lemma}

\begin{proof}
For each $n\geq1$, let $\lambda_n$ be the pushforward of
$\mu^{(n+2)}_{x_0,\ldots,x_{n+1}}$ under the map
$$(x_0,\ldots,x_{n+1}) \longmapsto (x_1,\ldots,x_n,x_0\cdots x_{n+1}),$$
and denote the last coordinate by $t$.
Idempotence of $\mu$,
associativity and functoriality of the Morley product imply that,
for each $1\leq i\leq n$, the map
$$(x_0,\ldots,x_{n+1}) \longmapsto (x_0\cdots x_{i-1},x_i,x_{i+1}\cdots x_{n+1})$$
pushes $\mu^{(n+2)}$ forward to $\mu^{(3)}$. Consequently, setting
$$\tau:=((x,y,z)\mapsto(y,xyz))_*(\mu^{(3)}),$$
and again using idempotence of $\mu$, we have
$$\lambda_n|_{x_1,\ldots,x_n}=\mu^{(n)},\qquad \lambda_n|_t=\mu,\qquad \lambda_n|_{x_i,t}=\tau\quad(1\leq i\leq n).$$
Lemma~\ref{lem:finite-couplings} gives $\tau=\mu^{(2)}$.
Interchanging the two coordinates and using symmetry of
$\mu^{(2)}$, i.e. Fact~\ref{fact:functorial},
proves the first assertion.

For $b\in G$, put $\nu_b:=\mu*\delta_b*\mu$.
This measure is definable over $Mb$, so
Fact~\ref{fact:uniqueness} gives
$$b\in K_\mu \quad\Longleftrightarrow\quad \nu_b|Mb=\mu|Mb.$$
Every formula over $Mb$ in the variable $t$ has the form
$\varphi(t,b)$ for some $\varphi(t,y)\in\CL(M)$.
Thus the condition on the right is equivalent to
$$(\mu_x\otimes\mu_z)\bigl(\varphi(xbz,b)\bigr) =\mu_t\bigl(\varphi(t,b)\bigr) \qquad\text{for every }\varphi(t,y)\in\CL(M).$$
For each $\varphi$, both sides are continuous functions of
$\tp(b/M)$, by definability of $\mu_x\otimes\mu_z$ and of $\mu$.
Their equality therefore defines a closed subset of $S_G(M)$.
Intersecting these closed sets shows that $K_\mu$ is
$M$-type-definable.

By symmetry of $\mu^{(3)}$ (Fact~\ref{fact:functorial}), the first assertion also gives
$$((x,z,y)\mapsto(xyz,y))_* \bigl((\mu_x\otimes\mu_z)\otimes\mu_y\bigr) =\mu_t\otimes\mu_y.$$
Applying Lemma~\ref{lem:parameter-pushforward} with
$f((x,z),y)=xyz$, we obtain
$$\mu*\delta_b*\mu=\mu \qquad\text{whenever }\tp(b/M)\in\supp(\mu|M).$$
Thus every such $b$ belongs to $K_\mu$. Since $K_\mu$ is
$M$-type-definable, $\mu$ concentrates on $K_\mu$.

Finally, let $\psi(x)\in\CL(M)$ define a superset of $K_\mu$.
By definability of $\mu$, the function
$$S_G(M)\ni \tp(b/M)\longmapsto(\mu*\delta_b)(\psi)\in [0,1]$$
is well defined and continuous. The convolution formula
(cf.~\cite[Remark~2.9]{HRz}) gives
$$1=\mu(\psi)=(\mu*\mu)(\psi) =\int_{S_G(M)}(\mu\ast\delta_b)(\psi)\,d(\mu|M).$$
Hence $(\mu\ast\delta_b)(\psi)=1$ on $\supp(\mu|M)$.
Since $\psi$ was arbitrary, compactness implies that
$\mu*\delta_b$ concentrates on $K_\mu$ whenever
$\tp(b/M)\in\supp(\mu|M)$.
\end{proof}

\begin{lemma}\label{lem:quotient-independence}
Let $h:G\to H$ be a definable homomorphism, let 
$\mu\in\mathfrak M_G(\FC)$ be fim and idempotent, and suppose that $\nu=h_*\mu$ is generically
transitive. Then
$$((x,y)\mapsto(xy,h(y)))_*(\mu^{(2)})=\mu\otimes\nu.$$
\end{lemma}
\begin{proof}
For $n\geq1$, consider the pushforward of $\mu^{(n+1)}$ under
the definable map
$$(x_0,\ldots,x_n)\longmapsto(P,u_1,\ldots,u_n),$$
where $P=x_0x_1\cdots x_n$, $u_i=h(x_i\cdots x_n)$ for $1\leq i\leq n$.
By Fact~\ref{fact:gt}, the map $(a,b)\mapsto(ab,b)$ preserves
$\nu^{(2)}$. 
For $1\leq j<n$, define
$$T_j:H^n\longrightarrow H^n,\qquad T_j(y_1,\ldots,y_n) =(y_1,\ldots,y_{j-1},y_jy_{j+1},y_{j+1},\ldots,y_n).$$
By Fact~\ref{fact:functorial} and the preceding identity,
$(T_j)_*(\nu^{(n)})=\nu^{(n)}$. Moreover,
$$(T_1\circ\cdots\circ T_{n-1})(y_1,\ldots,y_n) =(y_1\cdots y_n,y_2\cdots y_n,\ldots,y_n).$$
Since
$(h(x_1),\ldots,h(x_n))_*(\mu^{(n+1)})=\nu^{(n)}$,
we obtain 
$$
\begin{aligned}
(u_1,\ldots,u_n)_*(\mu^{(n+1)}) &=(T_1\circ\cdots\circ T_{n-1})_* \bigl((h(x_1),\ldots,h(x_n))_*(\mu^{(n+1)})\bigr)\\
&=(T_1\circ\cdots\circ T_{n-1})_*(\nu^{(n)})\\
&=\bigl((T_1)_*\circ\cdots\circ(T_{n-1})_*\bigr)(\nu^{(n)})\\
&=\nu^{(n)}.
\end{aligned}
$$

For $1\leq i\leq n$, put $A_i:=x_0\cdots x_{i-1}$, $B_i:=x_i\cdots x_n$.
Associativity and functoriality of the Morley product
(Fact~\ref{fact:functorial}), together with idempotence of $\mu$,
give
$$(A_i,B_i)_*(\mu^{(n+1)})=\mu^{(2)}.$$
As $P=A_iB_i$ and $u_i=h(B_i)$, we obtain
$$(P,u_i)_*(\mu^{(n+1)}) =((x,y)\mapsto(xy,h(y)))_*(\mu^{(2)}) =:\tau \qquad(1\leq i\leq n).$$
Now, let $\lambda_n(y_1,\ldots,y_n,t)$ be the pushforward of
$\mu^{(n+1)}$ under
$$(x_0,\ldots,x_n)\longmapsto(u_1,\ldots,u_n,P),$$
and let $s:G\times H\to H\times G$ be given by $s(a,b)=(b,a)$.
The preceding identities and idempotence of $\mu$ give
$$\lambda_n|_{y_1,\ldots,y_n}=\nu^{(n)},\qquad \lambda_n|_t=\mu,\qquad \lambda_n|_{y_i,t}=s_*\tau \quad(1\leq i\leq n).$$
Since $\nu=h_*\mu$ is fim by Fact~\ref{fact:functorial},
Lemma~\ref{lem:finite-couplings}, applied with $\nu$ in place of
$\mu$ and $\mu$ in place of $\nu$, yields $s_*\tau=\nu\otimes\mu$.
Finally, commutation of $\mu$ and $\nu$
(Fact~\ref{fact:functorial}) gives
$\tau=(s^{-1})_*(\nu\otimes\mu)=\mu\otimes\nu$,
as required.
\end{proof}

\begin{lemma}\label{lem:central-stationarity}
Let $\mu, \nu\in\mathfrak{M}_G(\FC)$ be fim.
Suppose that $\nu$ concentrates on a type-definable subgroup $H$,
is invariant under conjugation by $H$, and satisfies
$\mu*\nu=\mu$. Then $\nu$ concentrates on $\stab_{\mathrm r}(\mu)$.
\end{lemma}

\begin{proof}
Choose a small model $M$ over which the measures are fim and
$H$ is type-definable. For $n\geq2$, let $T_n:G^{n+1}\to G$
be the definable map
$$T_n(z,y_1,\ldots,y_n)=zy_1\cdots y_n,$$
and let $\lambda_{n;\, y_1,\ldots,y_n,t}$ be the pushforward of
$\mu_z\otimes\nu^{(n)}_{y_1,\ldots,y_n}$ under
$$(z,y_1,\ldots,y_n)\longmapsto (y_1,\ldots,y_n,T_n(z,y_1,\ldots,y_n)).$$
By the definition of $\lambda_n$, $\lambda_n|_{y_1,\ldots,y_n}=\nu^{(n)}$.
Moreover, associativity of convolution
(cf.~\cite[Theorem~2.18]{InTheWild}) and $\mu*\nu=\mu$ give
$$\lambda_n|_t =(T_n)_*\bigl(\mu_z\otimes\nu^{(n)}_{y_1,\ldots,y_n}\bigr) =\mu*\nu^{*n} =\mu.$$
Set $\tau:=((z,y)\mapsto(zy,y))_*(\mu\otimes\nu)$
and fix $1\leq i\leq n$. For every $b\in H$, conjugation invariance
of $\nu$ gives $\delta_b*\nu=\nu*\delta_b$. Hence
$$\mu*\nu^{*(i-1)}*\delta_b*\nu^{*(n-i)} =\mu*\nu^{*(n-1)}*\delta_b =\mu*\delta_b,$$
(where $\nu^{*0}:=\delta_1$).

To compute $\lambda_n|_{t,y_i}$, fix a formula
$\varphi(t,y;c)$ and a small model $N\supseteq Mc$.
By symmetry of the Morley powers of $\nu$ and associativity
(Fact~\ref{fact:functorial}),
$$\mu_z\otimes\nu^{(n)}_{y_1,\ldots,y_n} = \bigl( \mu_z\otimes \nu^{(n-1)}_{y_1,\ldots,y_{i-1},y_{i+1},\ldots,y_n} \bigr)\otimes\nu_{y_i}.$$
The definitions of pushforward, Morley product, and convolution
therefore give
$$
\begin{aligned}
\lambda_n(\varphi(t,y_i;c)) &= \int_{S_G(N)} \bigl(\mu*\nu^{*(i-1)}*\delta_b*\nu^{*(n-i)}\bigr) (\varphi(t,b;c))\,d(\nu|N)(b)\\
&= \int_{S_G(N)} (\mu*\delta_b)(\varphi(t,b;c))\,d(\nu|N)(b)\\
&=\tau(\varphi(t,y;c)).
\end{aligned}
$$
Here the second equality follows from the preceding convolution
identity, since $\nu$ concentrates on $H$. Thus
$\lambda_n|_{t,y_i}=\tau$ for every $1\leq i\leq n$.

Let $s:G^2\to G^2$ be given by $s(a,b)=(b,a)$. We have proved
$$\lambda_n|_{y_1,\ldots,y_n}=\nu^{(n)},\qquad \lambda_n|_t=\mu,\qquad \lambda_n|_{y_i,t}=s_*\tau \quad(1\leq i\leq n).$$
Since $\nu$ is fim, Lemma~\ref{lem:finite-couplings}, applied
with $\nu$ in place of $\mu$ and $\mu$ in place of $\nu$, yields
$s_*\tau=\nu\otimes\mu$.
Commutation of $\mu$ and $\nu$
(Fact~\ref{fact:functorial}) now gives
$\tau=(s^{-1})_*(\nu\otimes\mu)=\mu\otimes\nu$.

Applying Lemma~\ref{lem:parameter-pushforward} to
$f(z,b)=zb$, we obtain
$$\mu*\delta_b=\mu \qquad\text{whenever }\tp(b/M)\in\supp(\nu|M).$$
Since $\stab_{\mathrm r}(\mu)$ is $M$-type-definable,
$\nu$ concentrates on $\stab_{\mathrm r}(\mu)$.
\end{proof}

Taking $\nu=\mu$ in Lemma~\ref{lem:central-stationarity} shows
that every global fim idempotent measure $\mu$ which concentrates
on a type-definable subgroup $H$ and is invariant under
conjugation by $H$ is generically transitive.

\begin{lemma}\label{lem:relative-cancellation}
Let $h:G\to H$ be a definable surjective homomorphism and let
$\mu\in\mathfrak M_G(\FC)$ be fim and idempotent
such that $\nu=h_*\mu$ is generically
transitive. If
$$a\in C_G(\ker h),\qquad \mu*\delta_a*\mu=\mu,$$
then $\mu*\delta_a=\mu$.
\end{lemma}

\begin{proof}
Put $R=\stab_{\mathrm r}(\nu)$ and $L=h^{-1}(R)$.
Applying
Fact~\ref{fact:gt}(iv) gives $h(a)\in R$, hence $a\in L$.
Define maps $F_-,F_+\colon H\to G$ by
$$F_-(h(y))=y^{-1}ay, \qquad F_+(h(y))=yay^{-1}.$$
The maps $F_-$ and $F_+$
are well defined and definable because $a$ centralizes the normal
subgroup $\ker h$. Set $\eta=(F_-)_*\nu$.
This measure is fim and concentrates on $L$. For $b\in L$,
$$b^{-1}F_-(h(y))b=F_-(h(y)h(b)).$$
Right invariance of $\nu$ under $R$ therefore implies that $\eta$ is
invariant under conjugation by $L$.

Apply Lemma~\ref{lem:quotient-independence} to $h$. The map
$(t,s)\mapsto tF_-(s)$ sends $(xy,h(y))$ to $xay$, so pushing
forward gives
$$\mu*\eta=\mu*\delta_a*\mu=\mu.$$
Lemma~\ref{lem:central-stationarity} implies that $\eta$ concentrates
on $S=\stab_{\mathrm r}(\mu)$. Inversion preserves $\nu$, and
$F_+=F_-\circ\inv$, so $(F_+)_*\nu=\eta$ as well.
Put $A:=\{y\in G:yay^{-1}\in S\}
   =(F_+\circ h)^{-1}(S)$.
Since $(F_+\circ h)_*\mu=(F_+)_*\nu=\eta$
and $\eta$ concentrates on $S$, the measure $\mu$ concentrates
on $A$. For every $y\in A$, we have
$$\mu*\delta_{ya} =\mu*\delta_{yay^{-1}}*\delta_y =\mu*\delta_y.$$

Fix a formula $\varphi(t;c)$ and choose a small model $B$
containing $a,c$ such that $\mu$ is definable over $B$ and
$A$ is $B$-type-definable. By idempotence of $\mu$,
$$
\begin{aligned}
(\mu*\delta_a)(\varphi(t;c)) &=((\mu*\mu)*\delta_a)(\varphi(t;c))\\
&=(\mu*\mu)(\varphi(ta;c))\\
&=(\mu_x\otimes\mu_y)(\varphi((xy)a;c))\\
&=\int_{S_G(B)} \mu_x(\varphi(x(ya);c))\,d(\mu|B)(y)\\
&=\int_{S_G(B)} (\mu*\delta_{ya})(\varphi(t;c))\,d(\mu|B)(y)\\
&=\int_{S_G(B)} (\mu*\delta_y)(\varphi(t;c))\,d(\mu|B)(y)\\
&=(\mu*\mu)(\varphi(t;c))\\
&=\mu(\varphi(t;c)),
\end{aligned}
$$
where the sixth equality follows from $\mu$ concentrating on $A$.
As $\varphi(t;c)$ was arbitrary, $\mu*\delta_a=\mu$.
\end{proof}

\begin{theorem}\label{thm:commuting-kernels}
Let $N_1,N_2\trianglelefteq G$ be definable normal subgroups with
$[N_1,N_2]=1$, and let $h_i:G\to G/N_i$ be the quotient maps.
If 
$\mu\in\mathfrak M_G(\FC)$ is fim and idempotent
and both $(h_1)_*\mu$ and
$(h_2)_*\mu$ are generically transitive, then $\mu$ is generically
transitive.
\end{theorem}

\begin{proof}
Choose a small model $M$ over which $N_i$'s and $h_i$'s are definable,
and such that $\mu$ is fim over $M$. 
Consider
$$K:=\{g\in G:\mu*\delta_g*\mu=\mu\},\qquad \nu_i:=(h_i)_*\mu,\qquad R_i:=\stab_{\mathrm r}(\nu_i).$$
For $g\in K$, applying $(h_i)_*$ to $\mu*\delta_g*\mu=\mu$ gives
$\nu_i*\delta_{h_i(g)}*\nu_i=\nu_i$ for $i=1,2$.
Fact~\ref{fact:gt}(iv) therefore implies $h_i(g)\in R_i$.
Suppose that $g_1,g_2\in K$ and $h_1(g_1)=h_1(g_2)$.
Set
$$\lambda=\mu*\delta_{g_1},\qquad a=g_1^{-1}g_2\in N_1.$$
The measure $\lambda$ is fim and idempotent, and
$(h_2)_*\lambda=\nu_2$. Moreover,
$$\lambda*\delta_a*\lambda =\mu*\delta_{g_2}*\mu*\delta_{g_1} =\mu*\delta_{g_1}=\lambda.$$
Since $a$ centralizes $N_2$, Lemma~\ref{lem:relative-cancellation},
applied to $\lambda$ and $h_2$, yields
\begin{equation}
\label{eq:fibre-translates} \mu*\delta_{g_2}=\mu*\delta_{g_1}.
\end{equation}

Fix $b\models p\in\supp(\mu|M)$ and put
$\nu=\mu*\delta_b$. By Lemma~\ref{lem:sandwich-locus}, both
$\mu$ and $\nu$ concentrate on $K$, and $b\in K$. Hence
$$(h_1)_*\nu=\nu_1*\delta_{h_1(b)}=\nu_1=(h_1)_*\mu.$$
For a formula $\varphi(t;c)$, choose a model $N\supseteq Mbc$.
Let $S_K(N)$ be the closed subspace of $S_G(N)$ consisting of types
concentrating on $K$. The function
$$F:S_K(N)\longrightarrow[0,1],\qquad \tp(g/N)\longmapsto\mu_t(\varphi(tg;c))$$
is continuous and is constant on the fibres of
$\hat{h}_{1}:S_K(N)\longrightarrow S_{G/N_1}(N)$, given by
$\tp(g/N)\longmapsto\tp(h_1(g)/N)$.
Indeed, if $h_1(g_1)\equiv_N h_1(g_2)$, an automorphism over $N$ sends
$g_2$ to $g_2'$ with $h_1(g_2')=h_1(g_1)$. Invariance of $\mu$ and
\eqref{eq:fibre-translates} give $F(\tp(g_2/N))=F(\tp(g_1/N))$.

Since $S_K(N)$ is compact and $F$ is constant on the fibres of
$\hat{h}_1$, it induces a continuous function on the closed subset
$\hat{h}_1(S_K(N))$ of $S_{G/N_1}(N)$. By the Tietze extension
theorem, this function extends to a continuous function
$\hat{F}_1:S_{G/N_1}(N)\longrightarrow[0,1]$ such that $F=\hat{F}_1\circ\hat{h}_1$.
We compute
$$
\begin{aligned}
(\mu*\nu)(\varphi(t;c)) &=\int_{S_K(N)}F\,d(\nu|_N)\\
&=\int_{S_{G/N_1}(N)} \hat{F}_1\,d\bigl(((h_1)_*\nu)|_N\bigr)\\
&=\int_{S_{G/N_1}(N)} \hat{F}_1\,d\bigl(((h_1)_*\mu)|_N\bigr)\\
&=\int_{S_K(N)}F\,d(\mu|_N)\\
&=(\mu*\mu)(\varphi(t;c))\\
&=\mu(\varphi(t;c)).
\end{aligned}
$$
As $\varphi(t;c)$ was arbitrary, $\mu*\nu=\mu$. On the other hand,
associativity of convolution and idempotence of $\mu$ give
$$\mu*\nu =\mu*(\mu*\delta_b) =(\mu*\mu)*\delta_b =\mu*\delta_b =\nu.$$
Hence $\nu=\mu$. Every realization of every type in
$\supp(\mu|M)$ therefore belongs to $\stab_{\mathrm r}(\mu)$.
This proves generic transitivity.
\end{proof}

\begin{cor}\label{cor:centralizer-transfer}
Let $N\trianglelefteq G$ be definable. If a fim idempotent measure
$\mu\in\mathfrak{M}_G(\FC)$ 
has generically transitive pushforwards to both $G/N$ and
$G/C_G(N)$, then $\mu$ is generically transitive.
\end{cor}

\begin{proof}
The definable normal subgroups $N$ and $C_G(N)$ commute elementwise,
so Theorem~\ref{thm:commuting-kernels} applies.
\end{proof}

\section{Generic transitivity for virtually solvable groups}\label{sec:solvable}

\subsection{Solvable groups}
We first apply Corollary \ref{cor:centralizer-transfer} to abelian kernels.

\begin{cor}\label{cor:abelian-kernel-transfer}
Let $N\trianglelefteq G$ be definable and abelian, and let
$h:G\to G/N$ be the quotient map. If $\mu\in\mathfrak M_G(\FC)$
is fim and idempotent, and $h_*\mu$ is generically transitive,
then $\mu$ is generically transitive.
\end{cor}

\begin{proof}
Since $N\subseteq C_G(N)$, the quotient map $G\to G/C_G(N)$
factors through $h$. Generic transitivity is preserved by definable
homomorphisms, as follows from the shear characterization in
Fact~\ref{fact:gt}. Thus both quotient measures in
Corollary~\ref{cor:centralizer-transfer} are generically transitive.
\end{proof}

Recall that for a group $L$, its \emph{derived series} is defined recursively by
$$L^{(0)}=L,\qquad L^{(i+1)}=[L^{(i)},L^{(i)}] =\langle[g,h]:g,h\in L^{(i)}\rangle \quad(i\geq0).$$
The group $L$ is \emph{solvable} if $L^{(d)}=\{1\}$ for some
integer $d\geq0$. The least such $d$ is called the
\emph{derived length} of $L$. Thus the trivial group has derived
length zero, and every nontrivial abelian group has derived length one.

The \emph{derived words} are defined recursively by
$$\delta_0(x_1)=x_1,\qquad \delta_{i+1}(x_1,\ldots,x_{2^{i+1}}) =\bigl[\delta_i(x_1,\ldots,x_{2^i}), \delta_i(x_{2^i+1},\ldots,x_{2^{i+1}})\bigr] \quad(i\geq0).$$
For every $i\geq0$, the subgroup $L^{(i)}$ is generated by
the values of $\delta_i$ in $L$.

\begin{lemma}\label{lem:solvable-abelian-layer}
Let $N$ be a definable solvable group of derived length $d\geq 1$.
There is a definable characteristic abelian subgroup $A\leq N$
such that $N/A$ has derived length at most $d-1$.
\end{lemma}

\begin{proof}
Let $D$ be the definable set of values of $\delta_{d-1}$ in $N$,
and put
$$C=C_N(D),\qquad A=Z(C).$$
Both $C$ and $A$ are definable and characteristic in $N$.
The subgroup $N^{(d-1)}=\langle D\rangle$ is abelian, so it lies
in $C$. Every element of $C$ centralizes $D$, hence also
$N^{(d-1)}$. Therefore $N^{(d-1)}\subseteq Z(C)=A$, which proves
the assertion.
\end{proof}

\begin{theorem}\label{thm:solvable-kernel-transfer}
Let $N\trianglelefteq G$ be definable and solvable, and let
$h:G\to G/N$ be the quotient map. If $\mu\in\mathfrak M_G(\FC)$
is fim and idempotent and $h_*\mu$ is generically transitive,
then $\mu$ is generically transitive.
\end{theorem}

\begin{proof}
We induct on the derived length $d$ of $N$. The case $d=0$ is
immediate. For $d\geq 1$, choose $A$ as in
Lemma~\ref{lem:solvable-abelian-layer}. Since $A$ is characteristic
in $N$, it is normal in $G$. Let $\pi:G\to G/A$ and
$\bar h:G/A\to G/N$ be the quotient maps. By
Fact~\ref{fact:functorial}, $\pi_*\mu$ is fim and idempotent.
The kernel $N/A$ of $\bar h$ has derived length at most $d-1$,
and $\bar h_*(\pi_*\mu)=h_*\mu$. The induction hypothesis gives
generic transitivity of $\pi_*\mu$. Apply
Corollary~\ref{cor:abelian-kernel-transfer} to $\pi$.
\end{proof}

\begin{cor}[Solvable groups]\label{cor:solvable}
Every global fim idempotent Keisler measure on a definable solvable group
is generically transitive.
\end{cor}

\begin{proof}
Apply Theorem~\ref{thm:solvable-kernel-transfer} with $N=G$.
\end{proof}

\subsection{Virtually solvable groups}
We next pass to virtually solvable groups. Here virtual solvability
is an abstract group-theoretic condition: the solvable subgroup of
finite index need not be definable.

\begin{lemma}\label{lem:solvable-radical}
A definable virtually solvable group $L$ has a largest normal
solvable subgroup $K$. This subgroup is definable, characteristic,
and of finite index in $L$.
\end{lemma}

\begin{proof}
Let $S\leq L$ be a solvable subgroup of finite index $m=[L:S]$. Define
$$N_0:=\bigcap_{g\in L}gSg^{-1}.$$
This subgroup is the kernel of the action of $L$ on the finite
set of left cosets $L/S$ by left multiplication. Consequently, $N_0\trianglelefteq L$ and $[L:N_0]\leq m!$ and the subgroup $N_0$ is solvable.
The finite group
$L/N_0$ has a largest normal solvable subgroup 
(the product of its all finitely many
normal solvable subgroups is solvable).
Let $K$ be its preimage in $L$. Then $K$ is solvable, normal, and
of finite index. 
Every normal solvable subgroup of $L$ has solvable
image in $L/N_0$, so it is contained in $K$. In particular, $K$
is characteristic.

If $K=\{1\}$, it is definable. Otherwise, $K$ is the solvable
radical of $L$ and has finite derived length. By
\cite[Theorem~1.1(2)]{OuldHoucine2013}, it is definable in the
pure group $L$, and hence in the ambient structure.
\end{proof}

\begin{cor}\label{cor:finite-kernel-transfer}
Let $N\trianglelefteq G$ be finite, and let $h:G\to G/N$ be the
quotient map. If $\mu\in\mathfrak M_G(\FC)$ is fim and idempotent 
and $h_*\mu$ is generically transitive, then $\mu$ is generically transitive.
\end{cor}

\begin{proof}
Conjugation on $N$ embeds $G/C_G(N)$ into the finite group
$\Aut(N)$. The pushforward of $\mu$ to $G/C_G(N)$ is therefore
generically transitive by Fact~\ref{fact:finite}. Apply
Corollary~\ref{cor:centralizer-transfer}.
\end{proof}

\begin{theorem}
\label{thm:virtually-solvable-kernel-transfer}
Let $N\trianglelefteq G$ be definable and virtually solvable, and
let $h:G\to G/N$ be the quotient map. If
$\mu\in\mathfrak M_G(\FC)$ is fim and idempotent, and
$h_*\mu$ is generically transitive, then $\mu$ is generically
transitive.
\end{theorem}

\begin{proof}
By Lemma~\ref{lem:solvable-radical}, $N$ has a definable
characteristic solvable subgroup $K$ of finite index.
Let $\pi:G\to G/K$ and
$\bar h:G/K\to G/N$ be the quotient maps. 

The measure
$\pi_*\mu$ is fim and idempotent, the kernel $N/K$ of $\bar h$
is finite, and $\bar h_*(\pi_*\mu)=h_*\mu$ is generically transitive.
Corollary~\ref{cor:finite-kernel-transfer} gives generic
transitivity of $\pi_*\mu$. Theorem~\ref{thm:solvable-kernel-transfer},
applied to $\pi$, now gives generic transitivity of $\mu$.
\end{proof}

\begin{cor}[Virtually solvable groups]\label{cor:virtually-solvable}
Every global fim idempotent Keisler measure on a definable 
virtually solvable group
is generically transitive.
\end{cor}

\begin{proof}
Apply Theorem~\ref{thm:virtually-solvable-kernel-transfer} with $N=G$.
\end{proof}

\section{Generic transitivity along normal series}\label{sec:normal-series}
All definability and stable embeddedness assumptions are taken in the same
ambient theory $T^{\eq}$, over a fixed small set of parameters.

\begin{theorem}\label{thm:main}
Let $G$ be a definable group admitting a finite series
$$\{1\}=G_0\trianglelefteq G_1\trianglelefteq\cdots \trianglelefteq G_m=G$$
of definable subgroups normal in $G$. Suppose that, for each $i<m$, the
definable group $G_{i+1}/G_i$ satisfies at least one of the following
conditions:
\begin{enumerate}[(i)]
\item $G_{i+1}/G_i$ is virtually solvable;
\item $G_{i+1}/G_i$ is stably embedded and its induced structure is stable.
\end{enumerate}
Then every global fim idempotent measure on $G$ is generically transitive.
\end{theorem}

The proof uses the results for virtually solvable groups from the preceding
section and analogous results for stably embedded normal subgroups with
stable induced structure. We first establish the latter result. Throughout
the argument, we work in $T^{\eq}$. For a set $N$ definable over a small
base $B$, let $\mathcal N$ be its induced structure, as in
Section~\ref{sec:preliminaries}. 
As $\FC$ is $\kappa$-saturated and $|B|<\kappa$, we have that
$\mathcal N$ and $\mathcal N^{\eq}$ are $\kappa$-saturated.

The following lemma applies the standard coding of definable families
on stably embedded sets to the graphs of conjugation maps;
see \cite{Kamensky2009}.

\begin{lemma}\label{lem:stable-coding}
Let $N\trianglelefteq G$ be definable over $B$, and suppose that $N$ is
stably embedded with stable induced structure $\mathcal N$. There are a
group $D$ definable in $\mathcal N^{\eq}$ and an ambient definable
surjective homomorphism
$$r\colon G\longrightarrow D$$
with kernel $C_G(N)$. Thus $G/C_G(N)$ is definably isomorphic to $D$.
Moreover, for every $k\geqslant 1$, the ambient definable subsets of $D^k$
are exactly the subsets definable, with parameters, in
$\mathcal N^{\eq}$.
\end{lemma}

\begin{proof}
For $g\in G$, let
$$\Gamma_g=\{(u,v)\in N^2:v=gug^{-1}\}.$$
By the standard uniform coding of definable families on stably
embedded sets, using imaginaries (see \cite{Kamensky2009}), there
are a set $S$ definable in $\mathcal N^{\eq}$ and an ambient
definable map $r:G\to S$ such that, for all $g,h\in G$,
$$r(g)=r(h)\quad\Longleftrightarrow\quad\Gamma_g=\Gamma_h.$$
Put $D=r(G)$. Then $D$ is ambient
definable, and
$$
\begin{aligned}
r(g)=r(h) &\quad\Longleftrightarrow\quad gug^{-1}=huh^{-1}\text{ for every }u\in N\\
&\quad\Longleftrightarrow\quad h^{-1}g\in C_G(N).
\end{aligned}
$$
Since $C_G(N)$ is normal in $G$, the rule
$$r(g)r(h):=r(gh)$$
defines an ambient definable group operation on $D$.
Thus $r:G\to D$ is a surjective homomorphism with kernel $C_G(N)$.

We verify that $D$ and its group structure are definable in
$\mathcal N^{\eq}$. Choose a surjection $\pi:V\to S$ definable
in $\mathcal N^{\eq}$, where $V\subseteq N^m$ is definable in
$\mathcal N$. 
If $X\subseteq S^k$ is ambient definable,
then
$$(\pi^k)^{-1}(X)\subseteq V^k\subseteq N^{mk}$$
is ambient definable. 
By stable embeddedness, this inverse image is definable in
$\mathcal N$ with parameters from $N$. Replacing the imaginary
parameters defining $\pi$ by representatives in finite powers
of $N$, we may also regard $\pi$ as definable with parameters
from $N$. Since
$$X=\pi^k\bigl((\pi^k)^{-1}(X)\bigr),$$
it follows that $X$ is definable in $\mathcal N^{\eq}$ with
parameters from $N$.

Applying this observation to $D$ and to the graph of its
multiplication proves that $D$ is a group definable in
$\mathcal N^{\eq}$. 
For every $k\geqslant 1$, the same argument shows that each
ambient definable subset of $D^k$ is definable in
$\mathcal N^{\eq}$ with parameters from $N$.
Conversely, let $X\subseteq D^k$ be definable in
$\mathcal N^{\eq}$ with parameters. By the construction of
$\mathcal N^{\eq}$, after choosing representatives from $N$
for the imaginary parameters, the set $(\pi^k)^{-1}(X)\subseteq V^k$
is definable in $\mathcal N$ with parameters from $N$.
Since every basic relation of $\mathcal N$ is ambient
$B$-definable, this inverse image is ambient definable.
The map $\pi$ is also ambient definable, so $X=\pi^k\bigl((\pi^k)^{-1}(X)\bigr)$ is ambient definable. Thus the ambient definable subsets
of $D^k$ are exactly the subsets definable with parameters
in $\mathcal N^{\eq}$.
\end{proof}

\begin{lemma}\label{lem:stable-measure-transfer}
Let $N$ be a $B$-definable set which is stably embedded over $B$,
and suppose that its induced structure $\mathcal N$ over $B$
is stable.
Let $D$ be a group whose underlying set is contained in a finite
product of sorts of $\mathcal N^{\eq}$, viewed as ambient
definable sets in $\FC^{\eq}$. Suppose that $D$ and its
multiplication are definable over $B\cup c$ in $T^{\eq}$,
for some finite tuple $c$ from $N$.

For $\mu\in\mathfrak{M}_D(\FC)$, let $\bar\mu$
denote its restriction to formulas of $\mathcal N^{\eq}$.
If $\mu,\nu\in\mathfrak{M}_D(\FC)$ are definable, then
$$\overline{\mu\otimes\nu} =\bar\mu\otimes\bar\nu, \qquad \overline{\mu*\nu} =\bar\mu*\bar\nu.$$
Consequently, if $\mu\in\mathfrak{M}_D(\FC)$ is fim and idempotent, 
then $\mu$ is generically transitive in the ambient theory $T^{\eq}$.
\end{lemma}

\begin{proof}
Pulling back the domain of $D$ and the graph of its multiplication
to representatives in finite powers of $N$ gives $Bc$-definable
sets. By the definition of $\mathcal N$, these sets are definable
over $c$ in $\mathcal N$. Hence $D$ and its multiplication are
definable in $\mathcal N^{\eq}$.
The argument at the end of the proof of
Lemma~\ref{lem:stable-coding} shows that, for every $k\geqslant 1$,
the ambient definable subsets of $D^k$ are exactly the subsets
definable with parameters in $\mathcal N^{\eq}$.
Consequently, two ambient global measures on $D^k$ are equal
whenever their restrictions to $\mathcal N^{\eq}$ are equal.

Since $\mathcal N^{\eq}$ is stable, $\bar\mu$ and $\bar\nu$
are definable over small elementary submodels of
$\mathcal N^{\eq}$, by
\cite[Fact~2.2(1)(b), (3)(a)]{Artem_Kyle}.
Fix a formula $\varphi(x,y;d)$ of $\mathcal N^{\eq}$, with
$x,y$ ranging over $D$.
Choose a small elementary submodel
$A\preceq\mathcal N^{\eq}$ containing $c,d$, over which both
$\bar\mu$ and $\bar\nu$ are definable.
Choose a small ambient elementary model $M$ containing $B$
and representatives from $N$ of all elements of $A$, over
which both $\mu$ and $\nu$ are definable.

Every formula of $\mathcal N^{\eq}$ over $A$ has an ambient
translation over $M$. Thus the map
$$
\begin{aligned}
\rho\colon S_D(M) &\longrightarrow S_D^{\mathcal N^{\eq}}(A),\\
\tp(b/M) &\longmapsto \tp_{\mathcal N^{\eq}}(b/A)
\end{aligned}
$$
is well-defined and continuous.
For a formula $\psi(y;a)$ of $\mathcal N^{\eq}$ with
parameters $a$ from $A$, let $[\psi]_A$ and $[\psi]_M$
denote the clopen sets defined by $\psi$ in the respective
type spaces. Then
$$\rho^{-1}([\psi]_A)=[\psi]_M,$$
and therefore
$$
\begin{aligned}
\bigl(\rho_*(\nu|_M)\bigr)([\psi]_A) &=(\nu|_M)([\psi]_M)\\
&=\nu(\psi(y;a))\\
&=\bar\nu(\psi(y;a))\\
&=(\bar\nu|_A)([\psi]_A),
\end{aligned}
$$
so $\rho_*(\nu|_M)=\bar\nu|_A$.
Define
$$
\begin{aligned}
F\colon S_D^{\mathcal N^{\eq}}(A)&\longrightarrow[0,1],\\
\tp_{\mathcal N^{\eq}}(b/A) &\longmapsto\bar\mu_x(\varphi(x,b;d)).
\end{aligned}
$$
Definability of $\bar\mu$ over $A$ implies that $F$ is
well-defined and continuous. For every $b\in D$,
$$(F\circ\rho)(\tp(b/M)) =\bar\mu(\varphi(x,b;d)) =\mu(\varphi(x,b;d)).$$
The definitions of the two Morley products therefore give
$$
\begin{aligned}
\overline{\mu\otimes\nu}(\varphi(x,y;d)) &=(\mu\otimes\nu)(\varphi(x,y;d))\\
&=\int_{S_D(M)} F\circ\rho\,d(\nu|_M)\\
&=\int_{S_D^{\mathcal N^{\eq}}(A)} F\,d(\bar\nu|_A)\\
&=(\bar\mu\otimes\bar\nu)(\varphi(x,y;d)).
\end{aligned}
$$
Since $\varphi(x,y;d)$ was arbitrary, this proves $\overline{\mu\otimes\nu}=\bar\mu\otimes\bar\nu$.

Multiplication on $D$ is definable in both structures.
Hence, for every formula $\psi(z;e)$ of $\mathcal N^{\eq}$,
$$
\begin{aligned}
\overline{\mu*\nu}(\psi(z;e)) &=(\mu\otimes\nu)(\psi(xy;e))\\
&=(\bar\mu\otimes\bar\nu)(\psi(xy;e))\\
&=(\bar\mu*\bar\nu)(\psi(z;e)).
\end{aligned}
$$
This proves $\overline{\mu*\nu}=\bar\mu*\bar\nu$.
Finally, suppose that $\mu$ is fim and $\mu*\mu=\mu$.
The equality of convolutions gives $\bar\mu*\bar\mu=\overline{\mu*\mu}=\bar\mu$.
By Fact~\ref{fact:stable-imc}, there is a type-definable
subgroup $H\leqslant D$ in $\mathcal N^{\eq}$ such that
$\bar\mu$ concentrates on $H$ and is invariant under right
translation by every element of $H$.

Choose representatives from $N$ for a small set of parameters
over which $H$ is type-definable. The ambient translations of
its defining formulas show that $H$ is also ambient
type-definable.
If $X\subseteq D$ is ambient definable and contains $H$,
then $X$ is definable in $\mathcal N^{\eq}$, so
$\mu(X)=\bar\mu(X)=1$.
Thus $\mu$ concentrates on $H$.

For every ambient definable $X\subseteq D$ and every $b\in H$,
both $X$ and $Xb^{-1}$ are definable in $\mathcal N^{\eq}$.
Right invariance of $\bar\mu$ gives
$$(\mu*\delta_b)(X)=\mu(Xb^{-1})=\bar\mu(Xb^{-1})=\bar\mu(X)=\mu(X).$$
Hence $\mu*\delta_b=\mu$ for every $b\in H$, and therefore
$H\subseteq\stab_{\mathrm r}(\mu)$.
Since $\mu$ concentrates on $H$, it concentrates on
$\stab_{\mathrm r}(\mu)$.
Definition~\ref{def:gt} now implies that $\mu$ is generically
transitive in $T^{\eq}$.
\end{proof}

\begin{theorem}\label{thm:stable-kernel-transfer}
Let $N\trianglelefteq G$ be definable, stably embedded, and have stable
induced structure. Let $h\colon G\to G/N$ be the quotient map and let
$\mu\in\mathfrak{M}_G(\FC)$ be fim and idempotent. If $h_*\mu$ is
generically transitive, then $\mu$ is generically transitive.
\end{theorem}

\begin{proof}
By Lemma~\ref{lem:stable-coding}, the quotient $G/C_G(N)$ is definably
isomorphic to a group $D$ in the stable structure $\mathcal N^{\eq}$.
Let $r\colon G\to D$ be the resulting quotient homomorphism. Definable
pushforward preserves fim and commutes with convolution, so $r_*\mu$ is
fim and idempotent. 
Lemma~\ref{lem:stable-measure-transfer} shows that
it is generically transitive. Both $N$ and $C_G(N)$ are definable normal
subgroups of $G$, and they centralize one another. The generic
transitivity of the two quotient measures therefore implies that of
$\mu$, by Corollary~\ref{cor:centralizer-transfer}.
\end{proof}

\begin{proof}[Proof of Theorem~\ref{thm:main}]
We induct on the length $m$ of the series. The case $m=0$ is immediate.
Suppose $m>0$, put $N=G_1$, and let $h\colon G\to G/N$ be the quotient
map. The quotient has the definable normal series
$$\{1\}=G_1/N\trianglelefteq G_2/N\trianglelefteq\cdots \trianglelefteq G_m/N=G/N.$$
Its successive factors are canonically definably isomorphic to
$G_{i+1}/G_i$, for $1\leqslant i<m$. These isomorphisms are taken in $T^{\eq}$, so they preserve both alternative
hypotheses of the theorem.

Let $\mu\in\mathfrak{M}_G(\FC)$ be fim and idempotent. 
Its pushforward $h_*\mu$
is fim and idempotent, hence generically transitive by induction.
According to the hypothesis satisfied by $N$, apply
Theorem~\ref{thm:virtually-solvable-kernel-transfer} or
Theorem~\ref{thm:stable-kernel-transfer} to conclude that $\mu$ is
generically transitive.
\end{proof}

\appendix

\section{A counterexample for the automorphism case}
\label{sec:automorphism}

We follow the conventions of~\cite{GHK}.
Let $M\preccurlyeq\FC$, let $\bar m$ enumerate $M$, and let
$\bar x$ be the corresponding tuple of variables.
Write $S_{\bar m}(\FC)$ for the types extending
$\tp(\bar m/\emptyset)$ and
$\mathfrak M_{\bar m}(\FC)$ for the measures concentrated on this
space.

A measure is \emph{smooth over $M$} if its restriction to $M$
has a unique global extension.
For an infinite tuple of variables, a measure is fim over $M$
if every finite-coordinate marginal is fim over $M$ in the
sense of Definition~\ref{def:fim}; see
\cite[Definition~2.12]{GHK}.

For $\sigma\in\aut(\FC)$ and $\mu\in\mathfrak{M}_{\bar{x}}(\FC)$, 
the action is given by
$$(\sigma\cdot\mu)(\varphi(\bar x;\bar b)) =\mu(\varphi(\bar x;\sigma^{-1}(\bar b))), \qquad \stab(\mu) =\{\sigma\in\aut(\FC):\sigma\cdot\mu=\mu\}.$$
For an $M$-invariant measure $\mu$ and a larger monster model
$\FC'\succeq\FC$, write $\mu|_{\FC'}$ for its unique
$M$-invariant extension.

If $\mu\in\mathfrak M_{\bar m}(\FC)$ is Borel-definable over $M$
and $p=\tp(\sigma(\bar m)/\FC)$, where
$\sigma\in\aut(\FC')$, put
$$G_{\mu}^{\varphi(\bar x;\bar b)}(p) =(\mu|_{\FC'})(\varphi(\bar x;\sigma^{-1}(\bar b))).$$
This is a well-defined Borel function, and the convolution is given by
\begin{equation}
\label{eq:auto-convolution} (\mu*\nu)(\varphi(\bar x;\bar b)) =\int_{S_{\bar m}(\FC)} G_{\mu}^{\varphi(\bar x;\bar b)}(p)\,d\nu(p) \qquad (\nu\in\mathfrak M_{\bar m}(\FC));
\end{equation}
see \cite[Definition~4.18, Proposition~4.19, Definition~4.21 and
Remark~4.22]{GHK}.

For a partial type $\pi(\bar x;\bar y)$ over $\emptyset$ implying
$\bar x\equiv_\emptyset\bar y$, set
$$G_{\pi,\FC} =\{\sigma\in\aut(\FC): \FC\models\pi(\sigma(\bar m);\bar m)\}, \qquad [\pi(\bar x;\bar m)] =\{p\in S_{\bar m}(\FC):p\supseteq\pi(\bar x;\bar m)\}.$$
Conjecture~(A) of~\cite{GHK} asserts, in particular, that if
$\mu\in\mathfrak M_{\bar m}(\FC)$ is fim over $M$ and idempotent,
and $\stab(\mu)=G_{\pi,\FC}$, then
$\mu([\pi(\bar x;\bar m)])=1$.

\begin{theorem}\label{thm:auto}
Let $T=\operatorname{DLO}$ be the theory of dense linear orders
without endpoints in the language $\{<\}$, and let
$M=(\mathbb Q,<)\preccurlyeq\FC$.
There is a measure $\mu\in\mathfrak M_{\bar m}(\FC)$ which is
smooth over $M$ and fim over $M$, such that
$$\mu*\mu=\mu, \qquad \stab(\mu)=G_{\pi,\FC}, \qquad \mu([\pi(\bar x;\bar m)])=0$$
for a partial type $\pi$ over $\emptyset$.
Thus Conjecture~(A) of~\cite{GHK} fails in an NIP theory.
\end{theorem}

Index the enumeration by the underlying rational order:
$\bar m=(q)_{q\in\mathbb Q}$ and
$\bar x=(x_q)_{q\in\mathbb Q}$.
For $c\in\FC$, define
$$s(c)=\sup\{q\in\mathbb Q:q<c\} \in\mathbb R\cup\{-\infty,+\infty\}.$$
Moreover, we introduce
$$
F:\mathbb R\longrightarrow[0,1],\qquad F(t)= \begin{cases}
0,&t\leq0,\\
t,&0<t<1,\\
1,&t\geq1,
\end{cases}
$$
and additionally set $F(-\infty):=0$ and $F(+\infty):=1$.
Here, real numbers describe cuts in the fixed copy of $(\mathbb Q,<)$.

For $r\in(0,1)$, let $p_r^-(z)$ be the global unary type given by
$$c<z\quad\text{when }s(c)<r, \qquad z<c\quad\text{when }r\leq s(c).$$
Every finite set of lower bounds lies below every finite set of
upper bounds, so these conditions are consistent.
They determine a complete type by quantifier elimination.
Let $q_r(\bar x)$ be the complete global type determined by
$$\bigcup_{q\in\mathbb Q}p_r^-(x_q) \ \cup\ \{x_q<x_t:q<t\}.$$
Finite consistency follows by choosing finitely many points in
increasing order between the finitely many bounds; completeness
again follows from quantifier elimination.
In particular, $q_r\in S_{\bar m}(\FC)$.

Define
\begin{equation}
\label{eq:auto-measure} \mu(\varphi(\bar x;\bar b)) =\operatorname{Leb} \{r\in(0,1):\varphi(\bar x;\bar b)\in q_r\}.
\end{equation}
After quantifier elimination, the set on the right is a finite
Boolean combination of intervals, with endpoints among the cuts
$s(b_j)$.
Since each $q_r$ is a complete type,
\eqref{eq:auto-measure} defines a finitely additive probability
measure concentrated on $\tp(\bar m/\emptyset)$.

For $q,t\in\mathbb Q$ and $c\in\FC$, the definition gives
\begin{align}
\mu(x_q<c)&=F(s(c)), & \mu(x_q=c)&=0, \label{eq:auto-marginals}\\
\mu(x_q<x_t)&=1\quad(q<t), & \mu(x_q<c<x_t)&=0\quad(q<t). \label{eq:auto-collapse}
\end{align}

Fix $q_1,\ldots,q_n\in\mathbb Q$ with $q_1<\cdots<q_n$
and distinct parameters 
$b_1,\cdots,b_k\in\FC$ with
$b_1<\cdots<b_k$.
Put
$$
\begin{aligned}
C_0(\bar x)&:\quad x_{q_1}<\cdots<x_{q_n}<b_1,\\
C_j(\bar x)&:\quad b_j<x_{q_1}<\cdots<x_{q_n}<b_{j+1} \qquad(1\leq j<k),\\
C_k(\bar x)&:\quad b_k<x_{q_1}<\cdots<x_{q_n}.
\end{aligned}
$$
Additionally, in the case when there are no parameters $b_i$'s selected,
we set
$$C_0(\bar x):\quad x_{q_1}<\cdots<x_{q_n}.$$
Set $u_0=0$, $u_{k+1}=1$, and
$u_j=F(s(b_j))$ for $1\leq j\leq k$.
Then, for $0\leq j\leq k$,
\begin{equation}
\label{eq:auto-cells} \mu(C_j)=u_{j+1}-u_j, \qquad \mu\Big(\bigvee_{j=0}^{k}C_j\Big)=1.
\end{equation}
Indeed, outside the finite union of the zero-measure sets defined by
$x_{q_i}=b_\ell$ and $x_{q_1}<b_\ell<x_{q_n}$, all coordinates
belong to a single interval between consecutive parameters.
There, membership in $C_j$ is determined by the interval containing $x_{q_1}$. For $k\geq1$, \eqref{eq:auto-marginals} gives these intervals measures $u_{j+1}-u_j$, including the two unbounded intervals. When $k=0$, $C_0$ has measure $1$.
By quantifier elimination, every formula in these coordinates
over $\{b_1,\ldots,b_k\}$ is constant on each $C_j$.
Thus \eqref{eq:auto-cells} determines its measure.

\begin{lemma}\label{lem:auto-regularity}
The measure $\mu$ is smooth over $M$ and fim over $M$.
\end{lemma}

\begin{proof}
Fix rational numbers $q_1<\cdots<q_n$.
Let $N\succeq M$, and let $\rho$ be a measure over $N$ in the
variables $x_{q_1},\ldots,x_{q_n}$ which agrees with $\mu$ on every
formula over $M$ in these variables.
For $c\in N$ with $s(c)\in\mathbb R$, choose rationals
$a<s(c)<b$.
Then $a<c<b$, and hence
$$F(a)=\rho(x_{q_i}<a) \leq\rho(x_{q_i}<c) \leq\rho(x_{q_i}<b)=F(b).$$
Letting $a$ and $b$ tend to $s(c)$ gives
$$\rho(x_{q_i}<c)=F(s(c)).$$
Moreover,
$\{x_{q_i}=c\}\subseteq\{a<x_{q_i}<b\}$, so
$$\rho(x_{q_i}=c)\leq F(b)-F(a).$$
Thus $\rho(x_{q_i}=c)=0$.
At the two infinite cuts, the same conclusions follow from
$\rho(0<x_{q_i}<1)=1$.

For $i<j$, the formula $x_{q_i}<x_{q_j}$ has $\rho$-measure one (as $\mu$ and $\rho$ agree on $\CL(M)_{x_{q_1},\ldots,x_{q_n}}$-formulas).
On the clopen $[x_{q_i}<x_{q_j}]$,
$[x_{q_j}<c]\subseteq [x_{q_i}<c]$.
These two sets have equal measure, and
$\rho(x_{q_j}=c)=0$, so
$$\rho(x_{q_i}<c<x_{q_j})=0.$$
Consequently, the calculation in \eqref{eq:auto-cells} applies to
$\rho$ and determines its value on every formula over $N$.
This proves smoothness over $M$ of each finite-variable restriction.
Since every formula involves only finitely many variables,
$\mu$ is also smooth over $M$.

By \cite[Corollary~2.6]{HruPiSi13}, every finite-variable restriction is fim over $M$.
Hence $\mu$ is fim over $M$ by
\cite[Definition~2.12]{GHK}.
In particular, $\mu$ is definable over $M$.
\end{proof}

\begin{lemma}\label{lem:auto-absorption}
The measure $\mu$ satisfies $\mu*\mu=\mu$.
\end{lemma}

\begin{proof}
Let $p\in\supp(\mu)$, and choose a larger monster model
$\FC'\succeq\FC$ and $\sigma\in\aut(\FC')$ such that
$p=\tp(\sigma(\bar m)/\FC)$.
By \eqref{eq:auto-marginals}--\eqref{eq:auto-collapse}, for every
$c\in\FC$ the type $p$ contains
$$x_q\ne c\quad(q\in\mathbb Q), \qquad x_q<x_t,\quad \neg(x_q<c<x_t)\quad(q<t).$$
It follows that precisely one of the following alternatives holds:
$$c<\sigma(q)\text{ for every }q\in\mathbb Q, \qquad\text{or}\qquad \sigma(q)<c\text{ for every }q\in\mathbb Q.$$
Thus $\sigma^{-1}(c)$ lies either below every element of $M$
or above every element of $M$.

Write $\mu'=\mu|_{\FC'}$. For $c\in\FC$, the preceding dichotomy
gives either $\sigma^{-1}(c)<0$ and $c<\sigma(q)$,
or $1<\sigma^{-1}(c)$ and $\sigma(q)<c$.
Since $\mu'(0<x_q<1)=1$ and
$x_q<c\in p$ exactly when $\sigma(q)<c$, we obtain
$$\mu'(x_q<\sigma^{-1}(c)) =\mathbf1_{[x_q<c]}(p), \qquad \mu'(x_q=\sigma^{-1}(c))=0.$$
For every $q,t\in\mathbb Q$ with $q<t$, we have
$\mu'(x_q<x_t)=1$, and $\sigma$ preserves
the order of parameters.
Hence every atomic formula with parameters $\sigma^{-1}(\bar b)$
has $\mu'$-measure $1$ if the corresponding formula over $\bar b$
belongs to $p$, and $0$ otherwise.
These assignments hold simultaneously for any finite collection
of atomic formulas.
Quantifier elimination gives
$$G_{\mu}^{\varphi(\bar x;\bar b)}(p) =\mathbf1_{[\varphi(\bar x;\bar b)]}(p).$$
The regular Borel measure associated with $\mu$ gives
$\supp(\mu)$ measure one
\cite[Proposition~2.7(2)]{Artem_Kyle}.
Thus \eqref{eq:auto-convolution} yields
$$
\begin{aligned}
(\mu*\mu)(\varphi(\bar x;\bar b)) &=\int_{S_{\bar m}(\FC)} G_{\mu}^{\varphi(\bar x;\bar b)}(p)\,d\mu(p)\\
&=\int_{S_{\bar m}(\FC)} \mathbf1_{[\varphi(\bar x;\bar b)]}(p)\,d\mu(p)\\
&=\mu(\varphi(\bar x;\bar b)).
\end{aligned}
$$
Therefore $\mu*\mu=\mu$.
\end{proof}

\begin{proof}[Proof of Theorem~\ref{thm:auto}]
Lemmas~\ref{lem:auto-regularity} and~\ref{lem:auto-absorption}
show that $\mu$ is smooth over $M$, fim over $M$, and idempotent.
Since $\mu$ is definable over $M$,
\cite[Lemma~2.26]{GHK} gives a partial type
$\pi(\bar x;\bar y)$ over $\emptyset$, implying
$\bar x\equiv_\emptyset\bar y$, such that
$$\stab(\mu)=G_{\pi,\FC}.$$

Let $\sigma\in\stab(\mu)$.
For $t\in\{1/4,3/4\}$, since $\sigma\cdot\mu=\mu$, 
we have
$$F(s(\sigma(t))) =\mu(x_0<\sigma(t)) =(\sigma\cdot\mu)(x_0<\sigma(t)) =\mu(x_0<t)=t.$$
Since $F(u)=t$ with $0<t<1$ implies $u=t$, we obtain
$s(\sigma(t))=t$.
If $\sigma(1/4)\geq1/2$, the definition of $s$ would give
$s(\sigma(1/4))\geq1/2$, contradicting $s(\sigma(1/4))=1/4$.
Similarly, $\sigma(3/4)\leq1/2$ would imply
$s(\sigma(3/4))\leq1/2$, contradicting $s(\sigma(3/4))=3/4$.
Therefore,
$$\sigma(1/4)<1/2<\sigma(3/4).$$
We claim that
$$\pi(\bar x;\bar m)\vdash x_{1/4}<1/2<x_{3/4}.$$
Otherwise, saturation would give a tuple $\bar a$ in $\FC$
satisfying $\pi(\bar a;\bar m)$ and
$\neg(a_{1/4}<1/2<a_{3/4})$.
As $\bar a\equiv_\emptyset\bar m$, homogeneity would yield
$\sigma\in\aut(\FC)$ with $\sigma(\bar m)=\bar a$.
Then $\sigma\in G_{\pi,\FC}=\stab(\mu)$, contradicting the
preceding inequality.

It follows from \eqref{eq:auto-collapse} that
$$\mu([\pi(\bar x;\bar m)]) \leq \mu(x_{1/4}<1/2<x_{3/4})=0.$$
This contradicts the concentration conclusion of
Conjecture~(A) of~\cite{GHK}.
\end{proof}

\printbibliography

@article{Artem_Kyle,
	title        = {Definable convolution and idempotent {K}eisler measures},
	author       = {Chernikov, Artem and Gannon, Kyle},
	year         = 2022,
	journal      = {Israel J. Math.},
	volume       = 248,
	number       = 1,
	pages        = {271--314},
	doi          = {10.1007/s11856-022-2298-2},
	issn         = {0021-2172,1565-8511},
	fjournal     = {Israel Journal of Mathematics},
	mrclass      = {03C45 (03C60 28A05 37A05 37B05)},
	mrnumber     = 4429282,
	mrreviewer   = {Sergey\ V.\ Sudoplatov}
}

@article{ConGann20,
title = {Remarks on generic stability in independent theories},
journal = {Annals of Pure and Applied Logic},
volume = {171},
number = {2},
pages = {102736},
year = {2020},
doi = {10.1016/j.apal.2019.102736},
author = {Conant, Gabriel and Gannon, Kyle},
}

@phdthesis{KyleThesis,
	title        = {Approximation {T}heorems for {K}eisler {M}easures},
	author       = {Gannon, Kyle},
	year         = 2020,
	school       = {University of Notre Dame}
}

@article {Artem_Kyle2,
    AUTHOR = {Chernikov, Artem and Gannon, Kyle},
     TITLE = {Definable convolution and idempotent {K}eisler measures, {II}},
   JOURNAL = {Model Theory},
  FJOURNAL = {Model Theory},
    VOLUME = {2},
      YEAR = {2023},
    NUMBER = {2},
     PAGES = {185--232},
      ISSN = {2832-9058,2832-904X},
   MRCLASS = {03C45 (03C60 28D15 37B05 43A10)},
  MRNUMBER = {4662630},
       DOI = {10.2140/mt.2023.2.185},
}

@article {Anand_Udi2011,
    AUTHOR = {Hrushovski, Ehud and Pillay, Anand},
     TITLE = {On NIP and invariant measures},
   JOURNAL = {J. Eur. Math. Soc.},
    VOLUME = {13},
      YEAR = {2011},
    NUMBER = {4},
     PAGES = {1005--1061},
       DOI = {10.4171/JEMS/274}
}

@online{GHK,
    title={Convolution semigroups for automorphism dynamics}, 
    author={Kyle Gannon and Daniel Max Hoffmann and Krzysztof Krupiński},
    year={2025},
    eprinttype = {arxiv},
    eprint = {2507.23503},
}

@article{CGK,
    author  = {Chernikov, Artem and Gannon, Kyle and Krupi{\'n}ski, Krzysztof},
  title   = {Definable Convolution and Idempotent Keisler Measures III: Generic Stability, Generic Transitivity, and Revised Newelski's Conjecture},
  journal = {J. Lond. Math. Soc.},
  volume  = {114},
  number  = {1},
  year    = {2026},
  pages   = {e70639},
  doi     = {10.1112/jlms.70639}
}

@article {InTheWild,
    AUTHOR = {Conant, Gabriel and Gannon, Kyle and Hanson, James},
     TITLE = {Keisler measures in the wild},
   JOURNAL = {Model Theory},
    VOLUME = {2(1)},
      YEAR = {2023},
    PAGES = {1--67},
    doi = {10.2140/mt.2023.2.1},
}

@online{ConGanHan23,
    title={Generic stability, randomizations, and NIP formulas},
    author={Conant, Gabriel and Gannon, Kyle and Hanson, James E},
    eprinttype={arxiv},
    eprint={2308.01801},
    year={2023}
}

@article {Pym62,
    AUTHOR = {Pym, J. S.},
     TITLE = {Idempotent measures on semigroups},
   JOURNAL = {Pacific J. Math.},
  FJOURNAL = {Pacific Journal of Mathematics},
    VOLUME = {12},
      YEAR = {1962},
    doi = {10.2140/pjm.1962.12.685},
     PAGES = {685--698},
}

@article{HruPiSi13,
    title={Generically stable and smooth measures in NIP theories},
    author={Hrushovski, Ehud and Pillay, Anand and Simon, Pierre},
    journal={Transactions of the American Mathematical Society},
    volume={365},
    number={5},
    pages={2341--2366},
    year={2013},
    doi={10.1090/S0002-9947-2012-05626-1},
}

@article {kawada1940probability,
    AUTHOR = {Kawada, Yukiyosi and It\^o, Kiyosi},
     TITLE = {On the probability distribution on a compact group. {I}},
   JOURNAL = {Proc. Phys.-Math. Soc. Japan (3)},
  FJOURNAL = {Proceedings of the Physico-Mathematical Society of Japan. 3rd
              Series},
    VOLUME = {22},
      YEAR = {1940},
     PAGES = {977--998},
 doi = {10.11429/ppmsj1919.22.12_977},
}

@article {Wendel,
    AUTHOR = {Wendel, J. G.},
     TITLE = {Haar measure and the semigroup of measures on a compact group},
   JOURNAL = {Proc. Amer. Math. Soc.},
  FJOURNAL = {Proceedings of the American Mathematical Society},
    VOLUME = {5},
      YEAR = {1954},
     PAGES = {923--929},
 doi = {10.2307/2032557},
}

@article {Rudin,
    AUTHOR = {Rudin, Walter},
     TITLE = {Idempotent measures on {A}belian groups},
   JOURNAL = {Pacific J. Math.},
  FJOURNAL = {Pacific Journal of Mathematics},
    VOLUME = {9},
      YEAR = {1959},
    doi={10.2140/pjm.1959.9.195},
     PAGES = {195--209},
}

@article {Glicksberg2,
    AUTHOR = {Glicksberg, Irving},
     TITLE = {Convolution semigroups of measures},
   JOURNAL = {Pacific J. Math.},
  FJOURNAL = {Pacific Journal of Mathematics},
    VOLUME = {9},
      YEAR = {1959},
     PAGES = {51--67},
    doi = {10.2140/pjm.1959.9.51},
}

@online{HRz,
  author     = {Hoffmann, Daniel Max and Rzepecki, Tomasz},
  title      = {On Idempotent Measure Conjecture and Decomposition of Invariant Measures},
  year       = {2025},
  eprinttype = {arxiv},
  eprint     = {2511.22945}
}

@article{PillayStableEmbeddedness,
  author  = {Pillay, Anand},
  title   = {Stable embeddedness and {NIP}},
  journal = {The Journal of Symbolic Logic},
  volume  = {76},
  number  = {2},
  year    = {2011},
  pages   = {665--672},
  doi     = {10.2178/jsl/1305810769}
}

@online{CGH26,
  author      = {Conant, Gabriel and Gannon, Kyle and Hanson, James E.},
  title       = {Generically stable {Keisler} measures},
  year        = {2026},
  eprinttype  = {arxiv},
  eprint      = {2608.24605},
  eprintclass = {math.LO}
}

@article{OuldHoucine2013,
  author      = {Ould Houcine, Abderezak},
  title       = {A remark on the definability of the {Fitting} subgroup
                 and the soluble radical},
  journal     = {Mathematical Logic Quarterly},
  year        = {2013},
  volume      = {59},
  number      = {1--2},
  pages       = {62--65},
  doi         = {10.1002/malq.201200038},
  eprinttype  = {arxiv},
  eprint      = {1205.0573},
  eprintclass = {math.GR}
}

@article{Kamensky2009,
  author      = {Kamensky, Moshe},
  title       = {Definable groups of partial automorphisms},
  journal     = {Selecta Mathematica},
  year        = {2009},
  volume      = {15},
  number      = {2},
  pages       = {295--341},
  doi         = {10.1007/s00029-009-0530-0},
  eprinttype  = {arxiv},
  eprint      = {math/0607718},
  eprintclass = {math.LO}
}

@article{GannonYeDomination,
  author  = {Gannon, Kyle and Ye, Jinhe},
  title   = {An invitation to extension domination},
  journal = {Notre Dame Journal of Formal Logic},
  volume  = {64},
  number  = {3},
  year    = {2023},
  pages   = {253--280},
  doi     = {10.1215/00294527-2023-0005}
}

@article{GannonTransfer,
  author  = {Gannon, Kyle},
  title   = {On transfer maps and the {Morley} product
             in {NIP} theories},
  journal = {The Journal of Symbolic Logic},
  year    = {2026},
  pages   = {1--24},
  note    = {First View},
  doi     = {10.1017/jsl.2026.10201}
}
\end{document}